\def\titlep{Serre-Swan theorem for non-commutative  C$^{*}$-algebras. 
Revised edition\footnote{Original paper \cite{SS}. 
The essential mathematical statement is same as before.}}
\documentclass[11pt]{article}
\usepackage{graphicx,ifthen}
\font\germ=eufm10 at12pt
\def\goth#1{\hbox{\germ#1}}

\newcommand{\qed}{\hbox{\rule[-2pt]{3pt}{6pt}}}

\newtheorem{Thm}{Theorem}[section]

\newtheorem{ex}[Thm]{Example}
\newtheorem{defi}[Thm]{Definition}
\newtheorem{lem}[Thm]{Lemma}

\newtheorem{prop}[Thm]{Proposition}

\newcommand{\qedh}{\hfill\qed \\}

\newcommand{\vv}{\vspace{.3in}}
\newcommand{\ww}{\vv\noindent }

\def\kup{{\cal K}_{u}({\cal P})}
\def\pr{{\it Proof. }}
\def\exx{{\cal E}_{X}}
\def\gax{\Gamma_{X}}

\addtocounter{footnote}{1}
\def\sftt#1{
\setcounter{equation}{0}
\addtocounter{footnote}{1}
\section{#1}
}

\def\ssft#1{\subsection{#1}}

\def\cal#1{\mathcal #1}

\begin{document}
%
%
%
\def\autherp{Katsunori Kawamura}
\def\emailp{e-mail: kawamura@kurims.kyoto-u.ac.jp.}
\def\addressp{College of Science and Engineering Ritsumeikan University,\\
1-1-1 Noji Higashi, Kusatsu, Shiga 525-8577,Japan
}
%
%
\setcounter{footnote}{0}
\pagestyle{plain}
\setcounter{page}{1}
\setcounter{section}{0}

\begin{center}
{\Large \titlep}

\ww
\autherp
\footnote{\emailp}

\noindent
{\it \addressp}
\quad \\
\end{center}

%
%
\begin{abstract}
{\small
We generalize the Serre-Swan theorem to non-commutative C$^{*}$-algebras.
For a Hilbert C$^{*}$-module $X$ over a C$^{*}$-algebra ${\cal A}$,
we introduce a hermitian vector bundle $\exx$ associated to $X$.
We show that there is a linear subspace $\Gamma_{X}$ of
the space of all holomorphic sections of ${\cal E}_{X}$
and a flat connection $D$ on ${\cal E}_{X}$
with the following properties:
(i) $\Gamma_{X}$ is a Hilbert ${\cal A}$-module with
the action of ${\cal A}$ defined by $D$,
(ii) the C$^{*}$-inner product of $\Gamma_{X}$ is induced by
the hermitian metric of ${\cal E}_{X}$, 
(iii) ${\cal E}_{X}$ is isomorphic to an associated bundle 
of an infinite dimensional Hopf bundle,
(iv) $\Gamma_{X}$ is isomorphic to $X$.\par
}
\end{abstract}


%
%
\sftt{Introduction}\label{section:first}
The Serre-Swan theorem \cite{Karoubi,Serre,Swan} is described as follows:
%
%
\begin{Thm}
\label{Thm:ssc}
Let $\Omega$ be a connected compact Hausdorff space and
let $C(\Omega)$ be the algebra of all complex-valued continuous
functions on $\Omega$.
Assume that $X$ is a module over $C(\Omega)$. 
Then $X$ is finitely generated projective iff 
there is a complex vector bundle $E$ over $\Omega$ 
such that $X$ is isomorphic onto the module of all continuous sections of $E$.
\end{Thm}
By Theorem \ref{Thm:ssc},
finitely generated projective modules over 
the commutative C$^{*}$-algebra $C(\Omega)$ and complex vector bundles 
over $\Omega$ are in one-to-one correspondence up to isomorphism. 
In non-commutative geometry \cite{Connesbook,VariGau},
a certain module over a non-commutative C$^{*}$-algebra ${\cal A}$
is treated as a non-commutative vector bundle over 
the non-commutative space ${\cal A}$, 
generalizing Theorem \ref{Thm:ssc} in a sense of {\it point-less} geometry.
Therefore both a non-commutative space and a non-commutative
vector bundle are invisible even if one desires to look hard.

On the other hand,  
for a unital generally non-commutative C$^{*}$-algebra ${\cal A}$,
the functional representation on a certain geometrical space 
is studied by \cite{CMP94}. We review it as follows.
%
%
\begin{defi}
\label{defi:ukb}
A triplet $({\cal P},p,B)$ is the uniform K\"{a}hler bundle associated 
with ${\cal A}$
if ${\cal P}$ ($={\rm Pure}{\cal A}$) 
is the set of all pure states of ${\cal A}$,
endowed with the $w^{*}$-uniformity, i.e. the uniformity 
which induces the  $w^{*}$-topology,
$B$ ($={\rm Spec}{\cal A}$) is the spectrum of ${\cal A}$, 
the set of all equivalence classes of irreducible representations of ${\cal A}$, 
and $p$ is the natural projection from ${\cal P}$ onto $B$
by the GNS representation.
\end{defi}
For each $b\in B$, the fiber 
${\cal P}_{b}\equiv p^{-1}(b)$ is a K\"{a}hler manifold 
(Appendix D in \cite{CMP94}).
Especially, if ${\cal A}$ is commutative, 
then ${\cal P}\cong B$ and it is a compact Hausdorff space. 
In this case,
each fiber of $({\cal P},p,B)$ is a $0$-dimensional K\"{a}hler manifold.
Define $C^{\infty}({\cal P})$ the set of 
all fiberwise-smooth complex-valued functions on ${\cal P}$.
The product $*$ on $C^{\infty}({\cal P})$ is defined by
%
%
\begin{equation}
\label{eqn:star1} 
l*m\equiv l\cdot
m+ \sqrt{-1}X_{m}l\quad (l,m\in C^{\infty}({\cal P}))
\end{equation}
where $X_{l}$ is the holomorphic part of the complex Hamiltonian vector field 
of $l$ with respect to the K\"{a}hler form on ${\cal P}$.
Then $C^{\infty}({\cal P})$ is a $^*$ algebra with the unit ${\bf 1}$ and 
the involution $^{*}$ by complex conjugation,
which is not associative in general. 
Define the subset $C^{\infty}_{u}({\cal P})$ of $C^{\infty}({\cal P})$
consisting of uniformly continuous functions on ${\cal P}$. 
%
%
\begin{Thm}\label{Thm:th22}
For a unital non-commutative C$^{*}$-algebra ${\cal A}$, 
the Gel'fand representation 
%
%
\begin{equation}
\label{eqn:gelfand}
f_{A}(\rho)\equiv \rho(A)\quad (A\in {\cal A}, \,\rho\in {\cal P}),
\end{equation}
gives an injective $^*$ homomorphism $f$ from ${\cal A}$ 
into $C^{\infty}({\cal P})$
where $C^{\infty}({\cal P})$ is endowed with the $*$-product
in (\ref{eqn:star1}).
The norm $\|\cdot\|$ on $f({\cal A})$ defined by
%
%
\begin{equation}\label{eqn:normtwo}
{\displaystyle
\|l\|\equiv\sup_{\rho\in {\cal P}}\left|\left(\bar{l}*l\right)(\rho)
\right|^{\frac{1}{2}}
\quad(l\in f({\cal A})),
}
\end{equation}
is a C$^{*}$-norm on the associative $^*$ subalgebra $f({\cal A})$. 

Furthermore $f({\cal A})$ is precisely the subset $\kup$ of $
C^{\infty}_{u}({\cal P})$ defined by
%
%
\begin{equation}\label{eqn:eqn221}
\kup\equiv\{l\in C^{\infty}_{u}({\cal P}): 
\bar{l}*l, l*\bar{l}\in C_{u}^{\infty}({\cal P}),\,D^{2}l=\bar{D}^{2}l=0\}
\end{equation}
where $D$, $\bar{D}$ are the holomorphic
and anti-holomorphic part, respectively,
of covariant derivative of K\"{a}hler metric 
defined on each fiber of ${\cal P}$.
In consequence, the following equivalence of C$^{*}$-algebras holds:
\[ {\cal A}\cong\kup.\]
\end{Thm}
%
%
\pr See Proposition 3.2 in \cite{CMP94}. \qedh

\noindent
By Theorem \ref{Thm:th22},
it seems that there exists a {\it geometry consisting of points} 
associated with not only a commutative C$^{*}$-algebra
but also a non-commutative one.
According to Theorem \ref{Thm:th22}, 
we introduce a representation of a Hilbert C$^{*}$-module 
as the sections of a vector bundle over ${\cal P}$.

A vector space $X$ is a {\it Hilbert C$^{*}$-module} \cite{Jen-Thom,Pas}
over a C$^{*}$-algebra ${\cal A}$
if $X$ is a right ${\cal A}$-module with an ${\cal A}$-valued inner product 
$\langle \cdot|\cdot \rangle$ which satisfies 
$\langle \eta|\xi a\rangle =\langle \eta|\xi\rangle a$
for each $\eta,\xi\in X$ and $a\in {\cal A}$,
and $X$ is complete with respect to the norm $\|\cdot\|$ defined by
$\|\xi\|\equiv\|\langle \xi|\xi\rangle \|^{1/2}$ for $\xi\in X$.
%
%
\begin{defi}\label{defi:ato12}
The triplet $(\exx,\Pi_{X},{\cal P})$ is the atomic bundle associated with
a Hilbert C$^{*}$-module $X$ over a unital C$^{*}$-algebra ${\cal A}$ 
if it is the fiber bundle with the base space ${\cal P}$ and 
the total space $\exx$:
\[\exx\equiv \bigcup_{\rho\in {\cal P}}{\cal E}_{X,\rho}\]
where $\Pi_{X}$ is the natural projection from $\exx$ onto ${\cal P}$,
and the fiber ${\cal E}_{X,\rho}$ for $\rho\in {\cal P}$ is the Hilbert space
defined as follows:
Define the quotient vector space
${\cal E}^{o}_{X,\rho}\equiv X/N_{\rho}$
where $N_{\rho}$ is the closed subspace of $X$ defined by 
$N_{\rho}\equiv\{\, \xi\in X:\rho(\|\xi\|^{2})=0\,\}$.
Define the inner product $\langle \cdot|\cdot\rangle_{\rho}$
on ${\cal E}_{X,\rho}^{o}$ by 
%
%
\begin{equation}
\label{eqn:inner}
\langle [\xi]_{\rho}| [\eta]_{\rho}\rangle _{\rho}
\equiv \rho(\langle \xi|\eta\rangle )
\qquad(\,[\xi]_{\rho},[\eta]_{\rho}\in{\cal E}^{o}_{X,\rho}\,)
\end{equation}
where $[\xi]_{\rho}\equiv \xi+N_{\rho}\in{\cal E}^{o}_{X,\rho}$ for $\xi\in X$.
Let ${\cal E}_{X,\rho}$ denote the completion of ${\cal E}^{o}_{X,\rho}$
by the norm $\|\cdot\|_{\rho}$ associated with 
$\langle \cdot|\cdot\rangle _{\rho}$.
\end{defi}

We show the property of ${\cal E}_{X}$.
Let ${\cal H}$ denote a complex Hilbert space 
with $1\leq \dim {\cal H}\leq \infty$.
A triplet $(S({\cal H}),\mu,{\cal P}({\cal H}))$ is
{\it the Hopf (fiber) bundle over ${\cal H}$}
if {\it the projective Hilbert space} ${\cal P}({\cal H})$
and {\it the Hilbert sphere} $S({\cal H})$ are defined by
%
%
\begin{equation}
\label{eqn:hopfbundle}
{\cal P}({\cal H})\equiv ({\cal H}\setminus \{0\})/{\bf C}^{\times},\quad
S({\cal H})\equiv \{z\in {\cal H}:\|z\|=1\}
\end{equation}
and the projection $\mu$ from $S({\cal H})$ onto ${\cal P}({\cal H})$ is 
defined by $\mu(z)\equiv [z]$ for $z\in S({\cal H})$.
%
%
\begin{Thm}\label{Thm:geometry}
For $b\in B$ $(={\rm Spec}{\cal A})$, 
let ${\cal H}_{b}$ be a representative of $b$,
${\cal E}_{X}^{b}\equiv \Pi_{X}^{-1}({\cal P}_{b})$
and $\Pi_{X}^{b}\equiv \Pi_{X}|_{{\cal E}_{X}^{b}}$.
Then $({\cal E}_{X}^{b},\Pi_{X}^{b},{\cal P}_{b})$
is a locally trivial vector bundle which is isomorphic to
the associated bundle of $(S({\cal H}_{b}),\mu,{\cal P}({\cal H}_{b}))$ 
by a certain Hilbert space $F_{X}^{b}$.
\end{Thm}
One of our aims is a geometric realization of a Hilbert C$^{*}$-module.
We illustrate the two-step fibration structure of 
the atomic bundle as follows:

\noindent
%
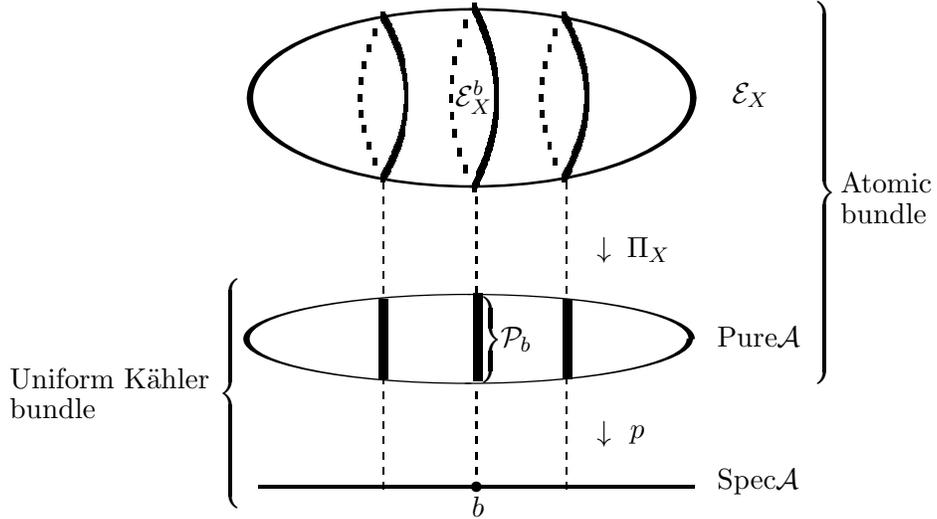
\def\tot{\scalebox{2}[.4]{
\linethickness{2mm}
\put(0,0){\circle{500}}
\linethickness{.2mm}}
}
\def\rot{
\qbezier[500](0,0)(100,100)(0,200)
\qbezier[10](0,0)(-100,100)(0,200)
}
\def\cuv{
\put(1000,0){\tot}
\put(800,-90){\scalebox{1}[.9]{\rot}}
\put(1000,-100){\rot}
\put(1200,-90){\scalebox{1}[.9]{\rot}}
}
\def\cuvtwo{\put(470,1000){\scalebox{3}[6]{\cuv}}}
%
%
\def\cuvthree{
\put(1050,0){\tot}
%
\put(1260,-90){\vbox{\hrule width 1pt height 10pt}}
\put(1060,-94){\vbox{\hrule width 1pt height 11pt}}
\put(850,-90){\vbox{\hrule width 1pt height 10pt}}
}
%
%
\def\base{
\put(300,1000){\scalebox{3}[3]{\cuvthree}}
\thicklines
\put(2050,15){\line(1,0){2900}}
\thinlines
%
\multiput(2880,0)(0,100){8}
{\line(0,1){45}}
\multiput(4100,0)(0,100){8}
{\line(0,1){45}}
\multiput(3500,0)(0,100){7}
{\line(0,1){45}}
%
\multiput(2880,1300)(0,100){8}
{\line(0,1){45}}
\multiput(4100,1300)(0,100){8}
{\line(0,1){45}}
\multiput(3500,1300)(0,100){7}
{\line(0,1){45}}
\put(400,670){Uniform K\"{a}hler}
\put(400,470){bundle}
\put(1800,1400){\rotatebox{-90}{$\underbrace{\hspace{1.2in}}$}}
\put(3660,950){${\cal P}_{b}$}
\put(3530,710){\rotatebox{90}{$\underbrace{\hspace{.45in}}$}}
\put(5100,950){${\rm Pure}{\cal A}$}
\put(5100,0){${\rm Spec}{\cal A}$}
\put(3470,-180){$b$}
\put(3500,20){\circle*{70}}
\put(4300,320){$\downarrow$}
\put(4520,340){$p$}
}
%
%
\def\module{
\put(1000,1000){\cuvtwo}
\put(6200,1970){${\cal E}_{X}$}
\put(4360,1950){${\cal E}_{X}^{b}$}
\put(6920,1370){Atomic}
\put(6920,1170){bundle}
\put(6750,100){\rotatebox{90}{$\underbrace{\hspace{2in}}$}}
\put(5300,940){$\downarrow$}
\put(5500,940){$\Pi_{X}$}
}
%
%
\setlength{\unitlength}{.0200mm}
\begin{picture}(3751,3551)(799,-200)
\thinlines
\put(0,0){\base}
\put(-1000,600){\module}
\end{picture}

Next, we reconstruct $X$ from ${\cal E}_{X}$.
Define the space of bounded sections
\[\Gamma(\exx)\equiv \{s:{\cal P}\to {\cal E}_{X}\,|\,
\Pi_{X}\circ s=id_{{\cal P}},\,\|s\|<\infty\}\]
where the norm $\|\cdot\|$ is defined by
%
%
\begin{equation}\label{eqn:norm}
\|s\|\equiv \sup_{\rho\in {\cal P}}\|s(\rho)\|_{\rho}.
\end{equation}
By standard operations, $\Gamma(\exx)$ is a complex linear space.
By Theorem \ref{Thm:geometry},
we can consider the differentiability of $s \in \Gamma(\exx)$ 
at each $B$-fiber in the sense of Fr\'echet differentiability 
of Hilbert manifolds.
Denote $\Gamma_{\infty}(\exx)$ the set of all
$B$-fiberwise smooth sections in $\Gamma(\exx)$.
Define the hermitian metric $H$ \cite{KobNo} 
on $\Gamma_{\infty}({\cal E}_{X})$ by 
%
%
\begin{equation}\label{eqn:hh1}
H_{\rho}(s,s^{'})\equiv\langle\, s(\rho)\,|\,s^{'}(\rho)\,\rangle _{\rho}
\quad(\rho\in{\cal P},\,s,s^{'}\in\Gamma_{\infty}(\exx)).
\end{equation}
By these preparations,
we state the following theorem which is a version of the Serre-Swan theorem
generalized to non-commutative C$^{*}$-algebras.
%
%
\begin{Thm}
\label{Thm:main00}
Let ${\cal A}$ be a unital C$^{*}$-algebra
with $({\cal P}, p, B)$ in Definition \ref{defi:ukb},
$f$ in (\ref{eqn:gelfand}) and $\kup$ in (\ref{eqn:eqn221}).
Let $X$ be a Hilbert ${\cal A}$-module with 
$(\exx,\Pi_{X},{\cal P})$ in Definition \ref{defi:ato12}
and $H$ in (\ref{eqn:hh1}).
Then the following holds:
\begin{enumerate}
\item
Let $X\times {\cal P}$ be the trivial bundle over ${\cal P}$ and
define the linear map $(P_{X})_{*}$ 
from $\Gamma(X\times {\cal P})$ to $\Gamma({\cal E}_{X})$ by 
$\{(P_{X})_{*}(s)\}(\rho)\equiv [s(\rho)]_{\rho}$
for $s\in \Gamma(X\times {\cal P}),\,\rho\in{\cal P}$.
Define the subspace $\Gamma_{X}$ of $\Gamma({\cal E}_{X})$ by
\[\gax\equiv (P_{X})_{*}(\Gamma_{const}(X\times {\cal P}))\]
where $\Gamma_{const}(X\times {\cal P})$ is the set of all constant 
sections of $X\times {\cal P}$.
Then any element in $\Gamma_{X}$ is holomorphic.
\item
There is a flat connection $D$ on ${\cal E}_{X}$
such that
$\gax$ is a Hilbert $\kup$-module with respect to the following right $*$-action 
%
%
\begin{equation}
\label{eqn:product}
s*l\equiv s\cdot l +\sqrt{-1}D_{X_{l}}s \qquad(\, (s,l)\in\gax\times\kup\,)
\end{equation}
and the C$^{*}$-inner product $H|_{\gax\times\gax}$.
\item
Under the identification $\kup$ with ${\cal A}$ by $f$,
the Hilbert ${\cal A}$-module $\gax$ is isomorphic to $X$.
\end{enumerate}
\end{Thm}

Here we summarize correspondences between geometry and algebra.

\noindent
\def\gel{
\put(300,1000){Gel'fand representation}
\put(0,0){
\begin{tabular}{lcc}
\hline
 & space & algebra \\\hline
 &       & $C(\Omega)$ \\
CG&$\Omega$& pointwise \\
&&product\\
\hline
NCG  & ${\cal P}\to B$
& $\kup$ \\
&& $*$-product\\
\hline
\end{tabular}
}}
\def\bundle{
\put(600,1000){Serre-Swan theorem}
\put(0,0){
\begin{tabular}{lcc}
\hline
& vector bundle& module \\
\hline
&  & $\Gamma(E)$ \\
CG& $E\to \Omega$ & pointwise \\
&& action \\
 \hline
NCG & $\exx\to {\cal P}$ & $\gax$ \\
&& $*$-action\\
\hline
\end{tabular}                              
}
}
%
%
\setlength{\unitlength}{.0200mm}
\begin{picture}(3751,2051)(100,-800)
\thinlines
\textsf{
\put(0,0){\gel}
\put(3000,0){\bundle}
}
\end{picture}

\noindent 
where we call respectively,
\textsf{CG} = commutative geometry as a geometry associated with commutative
C$^{*}$-algebras, and \textsf{NCG} = non-commutative geometry  
as a geometry associated with non-commutative C$^{*}$-algebras
according to \cite{Connes1}.
In this way, \textsf{NCG}'s are realized as visible geometries with points.

In $\S$ \ref{section:second}, 
we review the Hopf bundle and the uniform K\"{a}hler bundle.
In $\S$ \ref{subsection:csta}, we review \cite{CMP94} more closely. 
In $\S$ \ref{section:at1}, we show Theorem \ref{Thm:geometry}.
In $\S$ \ref{section:fourth},  we prove Theorem \ref{Thm:main00}.

%
%
\sftt{Hopf bundle and uniform K\"{a}hler bundle}
\label{section:second}

%
%
\ssft{The Hopf bundle and its associated bundle}
\label{subsection:hop}
We review the Hopf bundle and its associated bundle.
Let ${\bf S}\equiv (S({\cal H}),\mu,{\cal P}({\cal H}))$ be
the Hopf (fiber) bundle over a Hilbert space ${\cal H}$
in (\ref{eqn:hopfbundle}).
The space $S({\cal H})$ is a real submanifold of ${\cal H}$
in the relative topology.
We give ${\cal P}({\cal H})$ the quotient topology from
${\cal H}\setminus\{0\}\subset{\cal H}$ by the natural projection.
Then $\mu$ is continuous and open.

We define local trivial neighborhoods of the Hopf bundle according 
to Appendix C in \cite{CMP94}.
For $h\in S({\cal H})$, define
%
%
\begin{equation}
\label{eqn:coordinate}
\left\{
\begin{array}{c}
{\cal V}_{h}\equiv\{ [z]\in {\cal P}({\cal H}): \langle h|z\rangle \ne 0\},
\quad{\cal H}_{h}\equiv\{z\in {\cal H}:\langle h|z\rangle =0\},\\
\\
\beta_{h}:{\cal V}_{h}\to{\cal H}_{h};\quad 
\beta_{h}([z])\equiv \langle h|z\rangle^{-1}\cdot z-h\quad([z]\in{\cal V}_{h}).
\end{array}
\right.
\end{equation}
On the holomorphic tangent space $T_{\rho}{\cal P}({\cal H})$
at the local coordinate $({\cal V}_{h},\beta_{h},{\cal H}_{h})$
and $\beta_{h}(\rho)=z$,
we define 
the K\"{a}hler metric $g$ and 
the K\"{a}hler form $\omega$ on ${\cal P}({\cal H})$ by
\[g_{z}^{h}(\bar{v},u)
\equiv w_{z} \cdot\langle v|u\rangle
-w_{z}^{2}\cdot \langle v|z\rangle \langle z|u\rangle,
\quad  g_{z}^{h}(u,\bar{v})\equiv g_{z}^{h}(\bar{v},u),\]
\[\omega_{z}^{h}(\bar{v},u)\equiv 
\sqrt{-1}\{-w_{z} \cdot \langle v|u\rangle 
+w_{z}^{2}\cdot \langle v|z\rangle \langle z|u\rangle\},\quad
\omega_{z}^{h}(u,\bar{v})\equiv -\omega_{z}^{h}(\bar{v},u)\]
for $v,u\in {\cal H}_{h}$ where $w_{z}\equiv 1/(1+\|z\|^{2})$
and $\bar{x}\in {\cal H}_{h}^{*}$ means the dual vector of $x\in {\cal H}_{h}$.
Then ${\cal P}({\cal H})$ is a K\"{a}hler manifold with the holomorphic atlas
$\{({\cal V}_{h},\beta_{h},{\cal  H}_{h})\}_{h\in S({\cal H})}$.
For $l\in C^{\infty}({\cal P}({\cal H}))$,
define the {\it holomorphic Hamiltonian vector field} $X_{l}$ of $l$ 
by the equation
%
%
\begin{equation}
\label{eqn:hhm}
\omega_{\rho}((X_{l})_{\rho},\overline{Y}_{\rho})
=\bar{\partial}_{\rho}l(\overline{Y}_{\rho})\quad(
\overline{Y}_{\rho}\in \overline{T}_{\rho}{\cal P}({\cal H}),\,
\rho\in{\cal P}({\cal H}))
\end{equation}
where $\bar{\partial}$ is the anti-holomorphic differential operator
on $C^{\infty}({\cal P}({\cal H}))$ and $\overline{T}_{\rho}{\cal P}({\cal H})$
denotes the anti-holomorphic tangent space of ${\cal P}({\cal H})$
at $\rho\in {\cal P}({\cal H})$.

The family $\{{\cal V}_{h}\}_{h\in S({\cal H})}$ 
is a system of local trivial neighborhoods for ${\bf S}$ by 
the family $\{\psi_{h}\}_{h\in S({\cal H})}$ of maps defined by
$\psi_{h}:\mu^{-1}({\cal V}_{h})\to{\cal V}_{h}\times U(1)$;
%
%
\begin{equation}
\label{eqn:local}
\psi_{h}(z)\equiv (\,[z],\,\phi_{h}(z)\,),\quad
\phi_{h}(z)\equiv \langle z|h\rangle \cdot |\langle h|z\rangle |^{-1}.
\end{equation}
Furthermore we can verify that ${\bf S}$ is a principal $U(1)$-bundle.

Assume that $F$ is a complex vector space.
The fibration 
${\bf F}\equiv(\,S({\cal H})\times_{U(1)} \!F,\,\pi_{F},\, {\cal P}({\cal H})\,)$ is called {\it the associated bundle} of ${\bf S}$ by $F$
if $S({\cal H})\times_{U(1)} F$ is the set of all $U(1)$-orbits in
the product space $S({\cal H})\times F$ where the $U(1)$-action is defined by
\[(z,f)\cdot c\equiv \left(\bar{c}z,\,\bar{c}f\right)\qquad(\,c\in U(1),\,
(z,f)\in S({\cal H})\times F\,),\]
and the projection $\pi_{F}$ from $S({\cal H})\times_{U(1)}F$ onto 
${\cal P}({\cal H})$ is defined by $\pi_{F}([(x,f)])\equiv\mu(x)$
where we denote $[(x,f)]$ the element in $S({\cal H})\times_{U(1)} F$ 
containing $(x,f)$. 
The topology of $S({\cal H})\times_{U(1)} F$ is induced from
$S({\cal H})\times F$ by the natural projection.

For $h\in S({\cal H})$, 
the local trivialization $\psi_{F,h}$ of ${\bf F}$ at ${\cal V}_{h}$
is defined as the map
$\psi_{F,h}$ from $\pi_{F}^{-1}({\cal V}_{h})$ to ${\cal V}_{h}\times F$ by
%
%
\begin{eqnarray}
\label{eqnarry:ppp}
\psi_{F,h}([(z,f)])
\equiv&(\,\mu(z),\,\phi_{F,h}([(z,f)])\,),\quad
\phi_{F,h}([(z,f)])\equiv \phi_{h}(z)f.
\end{eqnarray}
The definition of $\psi_{F,h}$ is independent of the choice of $(z,f)$. 

%
%
\ssft{Connection}
\label{subsection:cona}
Let ${\bf F}=(\,S({\cal H})\times_{U(1)} \!F,\,\pi_{F},\, {\cal P}({\cal H})\,)$ 
be the associated bundle of the Hopf bundle ${\bf S}$ 
by $F$ in $\S$ \ref{subsection:hop}.
Let $\Gamma_{\infty}({\bf F})$ be the linear space of all smooth sections 
of ${\bf F}$.
A {\it connection} on ${\bf F}$ is a ${\bf C}$-bilinear map $D$ 
from ${\goth X}({\cal P}({\cal H})) \times\Gamma_{\infty}({\bf F})$ 
to $\Gamma_{\infty}({\bf F})$
which is $C^{\infty}({\cal P}({\cal H}))$-linear with respect to 
${\goth X}({\cal P}({\cal H}))$ and satisfies the 
Leibniz law with respect to $\Gamma_{\infty}({\bf F})$:
\[D_{Y}(s\cdot l)=\partial_{Y}l\cdot s+ l\cdot D_{Y}s
\quad(Y\in{\goth X}({\cal P}({\cal H})),\,
s\in \Gamma_{\infty}({\bf F}),\,l\in C^{\infty}({\cal P}({\cal H}))).\]
For $Y\in{\goth X}({\cal P}({\cal H}))$,
$h\in S({\cal H})$ and $\rho\in {\cal V}_{h}$, 
we denote $Y_{\rho}^{h}$  
the corresponding tangent vector at $\rho$ in a local chart.
Assume that a connection $D$ on ${\bf F}$ is written as
\[D=\partial+A.\]
According to the notation at the local chart,
we obtain families $\{A_{Y,\rho}^{h}:
Y\in{\goth X}({\cal P}({\cal H})),\,h\in S({\cal H}),\,\rho\in {\cal V}_{h}\}$ 
of linear maps  on $F$ such that
$\partial_{Y}|^{h}_{\rho}+A_{Y,\rho}^{h}
=(\partial_{Y}+A_{Y})_{\rho}^{h}=(\partial+A)_{Y,\rho}^{h}$.
Then we can verify that $D$ is a connection on ${\bf F}$ 
if and only if the following holds for each $h,h^{'}\in S({\cal H})$ 
with $<h|h^{'}>\ne 0$:
%
%
\begin{equation}
\label{eqn:cocycle}
A^{h^{'}}_{Y,\rho}
=-\frac{1}{2}\frac{\langle h|Y\rangle }{\langle h|z+h^{'}\rangle }
+A^{h}_{Y,\rho}\quad(\, \rho\in{\cal V}_{h^{'}}\cap {\cal V}_{h}\,)
\end{equation}
where $Y$ is a holomorphic tangent
vector of ${\cal P}({\cal H})$ at $\rho$ which is realized on
${\cal H}_{h^{'}}$ and $z=\beta_{h^{'}}(\rho)$.

A connection $D$ on ${\bf F}$ is {\it flat} if 
the {\it curvature} $R$ of ${\bf F}$ with respect to $D$ defined by
$R_{Y,Z}\equiv[D_{Y},\, D_{Z}]-D_{[Y,Z]}$,
$(Y,Z\in {\goth X}({\cal P}({\cal H}))$, vanishes.
%
%
\begin{prop}
\label{prop:conn2}
For $h\in S({\cal H})$ and
the chart $({\cal V}_{h}, \beta_{h},{\cal H}_{h})$ at $\rho\in{\cal P}({\cal H})$
in (\ref{eqn:coordinate}), 
we consider the trivializing neighborhood ${\cal V}_{h}$ for the Hopf bundle.
For $Y\in {\goth X}({\cal P}({\cal H}))$,
define the operator $D_{Y}$ on $\Gamma_{\infty}({\bf F})$ by 
\[(D_{Y}s)(\rho)\equiv (\partial_{Y}s)(\rho)+(A_{Y,\rho}s)(\rho)
\quad(\rho\in {\cal P}({\cal H}))\]
where $A_{Y,\rho}$ is defined as the family 
$\{A_{Y,\rho}^{h}:h\in S({\cal H}),\,\rho\in {\cal V}_{h}\}$
of linear operators on $F$ at $({\cal V}_{h}, \beta_{h},{\cal H}_{h})$, by
\[A_{Y,\rho}^{h}v\equiv 
-\frac{1}{2}\frac{\langle \beta_{h}(\rho)|Y_{\rho}^{h}\rangle }
{1+\|\beta_{h}(\rho)\|^{2}}\cdot v\quad (v\in F).\]
Then this defines a flat connection $D$ on ${\bf F}$.
\end{prop}
%
%
\pr
We can verify (\ref{eqn:cocycle}) for $\{A^{h}_{Y,\rho}\}$. 
Hence $D$ is a connection.
Furthermore it is straightforward to show that the curvature of $D$ vanishes.
\qedh

%
%
\ssft{Uniform K\"{a}hler bundle}\label{subsection:csta}
We show a geometric characterization of the set of all pure states 
and the spectrum of a C$^{*}$-algebra according to \cite{CMP94}.
%
%
\begin{defi}\label{defi:unik}
A triplet $(E,\mu,M)$ is called a uniform K\"{a}hler bundle 
if $E$ and $M$ are topological spaces and
$\mu$ is an open, continuous surjection from $E$ to $M$ such that
(i) the topology of $E$ is induced by a given uniformity,
(ii) each fiber $E_{m}\equiv \mu^{-1}(m)$ is a K\"{a}hler manifold.
\end{defi}
The local triviality of uniform K\"{a}hler bundle is not assumed. 
In general, the topological space $M$ is neither compact nor Hausdorff.

For uniform spaces, see Chapter 2 in \cite{Bourbaki}.  
Two uniform K\"{a}hler bundles $(E,\mu,M)$ and $(E^{'},\mu^{'},M^{'})$ are
{\it isomorphic} if there is a pair $(\beta,\phi)$ of a uniform  homeomorphism
$\beta$ from $E$ to $E^{'}$ and a homeomorphism
$\phi$ from $M$ to $M^{'}$, such that $\mu^{'}\circ\beta=\phi\circ \mu$
and any restriction $\beta |_{\mu^{-1}(m)} :\mu^{-1}(m)
\to (\mu^{'})^{-1}(\phi(m))$
is a holomorphic K\"{a}hler isometry for any $m\in M$.
We call $(\beta,\phi)$ a {\it uniform K\"{a}hler isomorphism} 
from $(E,\mu,M)$ to $(E^{'},\mu^{'},M^{'})$.

Let $({\cal H}_{b},\pi_{b})$ be an irreducible representation of ${\cal A}$
belonging to $b\in B$. 
Then $\rho\in {\cal P}_{b}$ corresponds $[x_{\rho}]\in{\cal P}({\cal H}_{b})$
where $\rho=\langle x_{\rho}| \pi_{b}(\cdot)x_{\rho}\rangle$.
Define the bijection $\tau^{b}$ 
from ${\cal P}_{b}$ onto ${\cal P}({\cal H}_{b})$ by
%
%
\begin{equation}\label{eqn:pureco}
\tau^{b}(\rho)\equiv[x_{\rho}]\quad(\rho \in {\cal P}_{b}).
\end{equation}
Then ${\cal P}_{b}$ has a K\"{a}hler manifold structure induced by $\tau^{b}$.
Furthermore the following holds.
%
%
\begin{Thm}
\label{Thm:thm21}
\begin{enumerate}
\item
For a unital C$^{*}$-algebra ${\cal A}$,
let $({\cal P},p,B)$ be as in Definition \ref{defi:ukb} and 
assume that $B$ is endowed with the Jacobson topology \cite{Ped}.
Then $({\cal P},p,B)$ is a uniform K\"{a}hler bundle.
\item
Let ${\cal A}_{i}$ be a C$^{*}$-algebra with the associated
uniform K\"{a}hler bundle $({\cal P}_{i}, p_{i}, B_{i})$ for $i=1,2$.
Then ${\cal A}_{1}$ and ${\cal A}_{2}$ are $^*$ isomorphic 
if and only if $( {\cal P}_{1}, p_{1},B_{1})$ 
and  $({\cal P}_{2},p_{2}, B_{2})$
are isomorphic as uniform K\"{a}hler bundle.
\end{enumerate}
\end{Thm}
%
%
\pr 
(i) See \cite{ACLM, CMP94}. 
(ii) See Corollary 3.3 in \cite{CMP94}.
\qedh
By Theorem \ref{Thm:thm21} (ii),
the uniform K\"{a}hler bundle $({\cal P},p, B)$ associated with
${\cal A}$ is uniquely determined up to uniform K\"{a}hler isomorphism.

By the above results, we obtain a fundamental correspondence between
algebra and geometry as follows:\\
\begin{tabular}{ccc}
&& \\
unital commutative C$^{*}$-algebra &
$\Leftrightarrow$ &
compact Hausdorff space \\
&&\\
$\bigcap$ & & $\bigcap$\\
&&\\
unital 
generally non-commutative  & 
$\Leftrightarrow$ & uniform K\"{a}hler bundle \\
C$^{*}$-algebra &&
associated with a C$^{*}$-algebra
\end{tabular}
\\
\\
The upper correspondence above is just the Gel'fand representation of
unital commutative C$^*$-algebras.
By these correspondences, we show the infinitesimal version of the 
Takesaki duality of Hamiltonian vector fields on a symplectic manifold \cite{Tak}.

%
%
\sftt{Proof of Theorem \ref{Thm:geometry}}\label{section:at1}
%
In this section,
we construct the typical fiber $F_{X}^{b}$ of ${\cal E}_{X}$
in Theorem \ref{Thm:geometry} and 
show the isomorphism among vector bundles.

In order to construct the typical fiber $F_{X}^{b}$ of ${\cal E}_{X}$,
we define the action $T=(t,\chi)$ of the group $G\equiv {\cal U}({\cal A})$
of all unitaries in ${\cal A}$ on $({\cal E}_{X},\Pi_{X},{\cal P})$ as follows:
The action $\chi$ of $G$ on the base space ${\cal P}$ is defined by
\[\chi_{u}(\rho)\equiv \rho \circ {\rm Ad}u^{*}
\quad(u\in G,\,\rho\in {\cal P}).\]
The action $t$ of $G$ on the total space ${\cal E}_{X}$ is defined by
\[t_{u}([\xi]_{\rho})\equiv
\left[\,\xi u^{*}\,\right]_{\chi_{u}(\rho)}\quad(u\in G,\, [\xi]_{\rho}
\in{\cal E}_{X,\rho}^{o}).\]
It is well-defined on the whole ${\cal E}_{X}$.
We see that $T=(t,\chi)$ is an action of $G$ on $(\exx,\Pi_{X},{\cal P})$
by bundle automorphism.
This action also preserves $B$-fibers $(\exx^{b},\Pi_{X}^{b},{\cal P}_{b})$
for each $b\in B$.

For $b\in B$, let $({\cal H},\pi)$ be a representative of $b$.
We identify ${\cal P}_{b}$ with ${\cal P}({\cal H})$ by
$\tau_{b}$ in (\ref{eqn:pureco}).
Furthermore we identify $\pi(u)$ with $u$ for each $u\in G$.
For the atomic bundle $(\exx^{b},\Pi_{X}^{b},{\cal P}_{b})$ and 
the Hopf bundle $(S({\cal H}),\mu_{b},{\cal P}_{b})$
in (\ref{eqn:hopfbundle}),
define their fiber product $\exx^{b}\times_{{\cal P}_{b}}S({\cal H})$ by 
\[\exx^{b}\times_{{\cal P}_{b}}S({\cal H})
= \{(x,h)\in\exx^{b}\times S({\cal H}):\Pi_{X}^{b}(x)=\mu_{b}(h)\}.\]
Thus the action $\sigma^{b}$ of $G$ 
on $\exx^{b}\times_{{\cal P}_{b}}S({\cal H})$ is defined by
\[\sigma_{u}^{b}(x,h)\equiv
\left(t_{u}(x), \pi_{b}(u)h\right)\qquad((x,h)\in
\exx^{b}\times_{{\cal P}_{b}}S({\cal H}),\,u\in G).\]
Define 
\[\mbox{$F_{X}^{b}$ the set of all orbits of $G$ 
in $\exx^{b}\times_{{\cal P}_{b}}S({\cal H})$}\]
and let ${\cal O}(x,h)\in F_{X}^{b}$ be the orbit of $G$ containing 
$(x,h)\in\exx^{b}\times_{{\cal P}_{b}}S({\cal H})$. 
We see that ${\cal O}(0,h)=\{(0,h^{'}):h^{'}\in S({\cal H})\}$.
We introduce the Hilbert space structure on $F_{X}^{b}$ as follows:
For $h\in S({\cal H})$,
define the sum and the scalar product on $F_{X}^{b}$ by
\[a{\cal O}(x,h)+b{\cal O}(y,h)\equiv {\cal O}(ax+by,h)\quad
(a,b\in {\bf C},\, x,y\in {\cal E}_{X}^{b}).\]
Then this operation is independent in the choice of $x,y$ and $h$.
For $h\in S({\cal H})$,
define the inner product $\langle \cdot|\cdot\rangle $
on the vector space $F_{X}^{b}$ by
\[\langle {\cal O}(x,h)|{\cal O}(y,h)\rangle \equiv \langle x|y\rangle _{\rho}
\quad(x,y\in {\cal E}_{X}^{b})\]
where $\rho=\mu_{b}(h)$.
Then $\langle {\cal O}(x,h)|{\cal O}(y,h)\rangle $
is independent in the choice of $x,y,\rho$ and $h$.
For $h_{0}\in S({\cal H})$ with $\mu_{b}(h_{0})=\rho$,
define the map $R_{\rho}$ from ${\cal E}_{X,\rho}$ to $F_{X}^{b}$ by 
$R_{\rho}(x)\equiv {\cal O}(x,h_{0})$ for $x\in {\cal E}_{X,\rho}$.
Then $R_{\rho}$ is a unitary from ${\cal E}_{X,\rho}$ to $F_{X}^{b}$
for each $\rho\in{\cal P}_{b}$.
In this way, $F_{X}^{b}$ is a Hilbert space.

We introduce the Hilbert bundle isomorphism in Theorem \ref{Thm:geometry}.
Let ${\bf F}_{X}^{b}\equiv(S({\cal H})
\times _{U(1)} F_{X}^{b},\,\pi_{F_{X}^{b}},
\,{\cal P}({\cal H})\,)$ be the associated bundle of
$(S({\cal H}),\mu_{b},{\cal P}({\cal H}))$ by $F_{X}^{b}$.
%
%
\begin{lem}
\label{lem:representation}
Any element of $S({\cal H})\times _{U(1)} F_{X}^{b}$
can be written as $[(h,{\cal O}(x,h))]$ where ${\cal O}(x,h)\in F_{X}^{b}$.
\end{lem}
%
%
\pr
By definition of the associated bundle in $\S$ \ref{subsection:hop},
an element of $S({\cal H})\times _{U(1)} F_{X}^{b}$ is  
the $U(1)$-orbit $[(h,{\cal O}(y,k))]$.
Because $({\cal H},\pi)$ is an irreducible representation of ${\cal A}$,
the action of $G$ on $S({\cal H})$ is transitive.
By this and definition of ${\cal O}(y,k)$,
there is $u\in G$ such that $h=uk$ and $(t_{u}^{b}(y),h)\in{\cal O}(y,k)$.
Denote $x\equiv t_{u}(y)$.
Then ${\cal O}(x,h)={\cal O}(y,k)$.
Hence $[(h,{\cal O}(y,k))]=[(h,{\cal O}(x,h))]$.
\qedh

\noindent
%
%
{\it Proof of Theorem \ref{Thm:geometry}.}
By Lemma \ref{lem:representation}, we shall denote
\[[h,x]\equiv[(h,{\cal O}(x,h) )]\in S({\cal H})\times _{U(1)} F_{X}^{b}
\qquad (h\in S({\cal H}),\,x\in\exx^{b}).\]
Define the map $\Phi^{b}$ from $\exx^{b}$ to 
$S({\cal H})\times_{U(1)}F_{X}^{b}$ by 
\[\Phi^{b}(x)\equiv[h_{x},x]\quad(x\in\exx^{b})\]
where $h_{x}\in \mu_{b}^{-1}(\Pi_{X}^{b}(x))$.
By definition of $F_{X}^{b}$, the map $\Phi^{b}$ is bijective.
We obtain a set-theoretical isomorphism $(\Phi^{b},\tau^{b})$ 
of fibrations between $(\exx^{b}$, 
$\Pi_{X}^{b}$, ${\cal P}_{b})$ and ${\bf F}_{X}^{b}$
such that any restriction
$\Phi^{b}|_{{\cal E}_{X,\rho}}$ of $\Phi^{b}$ at a fiber
${\cal E}_{X,\rho}$ is a unitary from ${\cal E}_{X,\rho}$ to
$\pi_{F_{X}^{b}}^{-1}(\rho)$ for $\rho\in {\cal P}_{b}$.
This unitary induces the Hilbert bundle isomorphism 
from $(\exx^{b},\Pi_{X}^{b},{\cal P}_{b})$ to ${\bf F}_{X}^{b}$.
\qedh

%
%
\sftt{Proof of Theorem \ref{Thm:main00}}\label{section:fourth}
Let us summarize our notations.
Let ${\cal A}$ be a unital C$^{*}$-algebra 
with the uniform K\"{a}hler bundle $({\cal P},p,B)$ and
let $X$ be a Hilbert C$^{*}$-module over ${\cal A}$
with the atomic bundle $\exx=(\exx,\Pi_{X},{\cal P})$.

Fix $b\in B$ and assume that $({\cal H},\pi)$ is a representative of $b$.
For the Hilbert space ${\cal H}$,
let $\{({\cal V}_{h},\beta_{h},{\cal H}_{h})\}_{h\in S({\cal H})}$
be as in (\ref{eqn:coordinate}).
For $\rho\in {\cal V}_{h}$,
define the vector $\Omega_{\rho}^{h}$ in ${\cal H}$ by
\[\Omega_{\rho}^{h}\equiv 
\{1+\|\beta_{h}(\rho)\|^{2}\}^{-1/2}\cdot \{\beta_{h}(\rho)+h\}.\]
Then $\rho=\langle \Omega_{\rho}^{h}|\pi(\cdot)\Omega_{\rho}^{h}\rangle$ and
$\langle h|\Omega_{\rho}^{h}\rangle \,> \,0$.
We prepare two lemmata to prove Theorem \ref{Thm:main00}. 
%
%
\begin{lem}
\label{lem:seceqn1}
For $s\in\Gamma({\cal E}_{X})$,
assume that there is a family $\{\xi_{\rho}\in X:\rho\in {\cal P}\}$ 
such that $s(\rho)=[\xi_{\rho}]_{\rho}\in {\cal E}_{X,\rho}$
for each $\rho\in {\cal P}$
and we identify ${\cal E}_{X}^{b}$ with $S({\cal H})\times_{U(1)}F_{X}^{b}$
by Theorem \ref{Thm:geometry}. 
Let $z=\beta_{h}(\rho)$ for $h\in S({\cal H})$
such that $\rho \in {\cal V}_{h}$.
Define $w_{z}\equiv 1/(1+\|z\|^{2})$ and
let $\phi_{F,h}$ be as in (\ref{eqnarry:ppp}) for $F=F_{X}^{b}$.
Then the following equations hold:
%
%
\begin{equation}
\label{eqn:local2}
\langle \,e\,|
\,\phi_{F,h}(s(\rho))\,\rangle 
=\sqrt{w_{z}} \cdot \langle \Omega_{\rho^{'}}^{h}
|\pi(\langle \xi^{'}|\xi_{\rho}\rangle )(z+h)\rangle,
\end{equation}
%
%
\begin{equation}
\label{eqn:deri22}
\partial_{Y}\phi_{F,h}(s(\rho))
={\cal O}(\,[
\partial_{Y}\hat{\xi}_{\rho}+\xi_{\rho}
(K_{Y,\rho}^{h}-2^{-1}w_{z}\langle z|Y\rangle )
]_{\rho},\, h\,)
\end{equation}
for $e={\cal O}([\xi^{'}]_{\rho^{'}},h)\in F_{X}^{b}$
where $K_{Y,\rho}^{h}\in {\cal A}$  is defined by
%
%
\begin{equation}\label{eqn:deri23}
\pi(K_{Y,\rho}^{h})(h+z)=Y
\end{equation}
and $[\partial_{Y}\hat{\xi}_{\rho}]_{\rho}\in {\cal E}_{X,\rho}$
is  defined by
$\langle \,[\eta]_{\rho}\,|\,[\partial_{Y}
\hat{\xi}_{\rho}]_{\rho}\,\rangle _{\rho}
\equiv \rho(\partial_{Y}\langle \eta|\xi_{\rho}\rangle )$
for $[\eta]_{\rho}\in{\cal E}_{X,\rho}$.
\end{lem}
%
%
\pr
By definition, we have that
$\phi_{F,h}(s(\rho))=c_{z,h}\cdot {\cal O}([\xi_{\rho}]_{\rho},z)$
where $c_{z,h}\equiv \langle z|h\rangle \cdot |\langle h|z\rangle |^{-1}$.
We have 
\[\langle \,e\,|\,\phi_{F,h}(s(\rho))\,\rangle 
=c_{z,h}\langle{\cal O}([\xi^{'}]_{\rho^{'}},h)|
{\cal O}([\xi_{\rho}]_{\rho},z_{\rho})\rangle.\]
Let $u\in G$ such that $\pi(u^{*})z=h=\Omega_{\rho^{'}}^{h}$.
Then 
${\cal O}([\xi_{\rho}]_{\rho},z)={\cal O}([\xi_{\rho}u]_{\rho^{'}},\pi(u^{*})z)$.
By this,
\[\langle{\cal O}([\xi^{'}]_{\rho^{'}},h)|
{\cal O}([\xi_{\rho}]_{\rho},z_{\rho})\rangle
=
\langle\Omega_{\rho^{'}}^{h}|
\pi_{b}(\langle\xi^{'}|\xi_{\rho}\rangle)\pi_{b}(u)\Omega_{\rho^{'}}^{h}\rangle\\
=\langle\Omega_{\rho^{'}}^{h}|
\pi_{b}(\langle\xi^{'}|\xi_{\rho}\rangle)z_{\rho}\rangle.\]
Because $z_{\rho}=c_{h,z}\Omega_{\rho}^{h}$,
(\ref{eqn:local2}) is verified.

By (\ref{eqn:local2}), we get
\[
\begin{array}{rl}
\langle \,e\,|\,\partial_{Y}
\phi_{F,h}(s(\rho))\,\rangle 
=& 
\sqrt{w_{z}}\cdot[ \langle \Omega_{\rho^{'}}^{h}|
\pi(\partial_{Y}\langle \xi^{'}|
\xi_{\rho}\rangle )(z+h)\rangle+
\langle \Omega_{\rho^{'}}^{h}|\pi(\langle \xi^{'}|
\xi_{\rho}\rangle )Y\rangle] \\
&-2^{-1}w_{z}^{3/2}\cdot \langle \Omega_{\rho^{'}}^{h}|\pi(\langle \xi^{'}|
\xi_{\rho}\rangle )(z+h)\rangle \langle z|Y\rangle.\\
\end{array}
\]
Hence we obtain (\ref{eqn:deri22}).
\qedh

For $\xi\in X$, define the section $s_{\xi}$ of $\exx$ by
$s_{\xi}(\rho)\equiv [\xi]_{\rho}$ for $\rho\in {\cal P}$.
Then $\|s_{\xi}\|= \|\xi\|$ for every $\xi\in X$.
Define the linear isometry $\Psi$ from $X$ into $\Gamma(\exx)$ by
\[\Psi(\xi)\equiv s_{\xi}\quad(\xi\in X).\]
%
%
\begin{lem}
\label{lem:mlem}
\begin{enumerate}
\item
For each $\xi\in X$, $\Psi(\xi)$ belongs to $\Gamma_{\infty}(\exx)$ and is holomorphic.
\item
According to Theorem \ref{Thm:geometry}, 
define the connection $D$ on ${\cal E}_{X}$
by the one in Proposition \ref{prop:conn2} at each fiber.
Let $*$ be as in (\ref{eqn:product}) with respect to $D$.
Then $\Psi(\xi)*f_{A}=\Psi(\xi\cdot A)$ for $\xi\in X$ and  $A\in {\cal A}$.
\end{enumerate}
\end{lem}
%
%
\pr
Let $\rho\in {\cal P}_{b}$ for $b\in B$.
Choose as a representative for $b$ an irreducible representation
$({\cal H}, \pi)$.
Fix $h\in S({\cal H})$ and, using the notations in (\ref{eqnarry:ppp}),
take the local trivialization $\psi_{F,h}$ of the Hopf bundle at 
$({\cal V}_{h},\beta_{h}, {\cal H}_{h})$ with $\rho \in {\cal V}_{h}$. 
Let $z\equiv\beta_{h}(\rho)\in {\cal H}_{h}$ and
$w_{z}\equiv 1/(1+\|z\|^{2})$.

\noindent
(i) Applying (\ref{eqn:deri22}) for $s=s_{\xi}$, we obtain
%
%
\begin{equation}\label{eqn:deri34}
\partial_{Y}\phi_{F,h}(s_{\xi}(\rho))
={\cal O}([\partial_{Y}\hat{\xi}+
\xi(K_{Y,\rho}^{h}-2^{-1}w_{z}\cdot 
\langle z|Y\rangle )]_{\rho}, \, h\,).
\end{equation}
Owing to (\ref{eqn:deri23}), the right-hand side of (\ref{eqn:deri34})
is smooth with respect to $z$. 
Hence $s_{\xi}$ is smooth at ${\cal P}_{b}$ for each $b\in B$. 
For $\rho_{0} \in{\cal P}_{b}$, we can choose $h_{0}\in S({\cal H})$ 
such that $\rho_{0}=\langle h_{0}|\pi(\cdot)h_{0}\rangle$.
Then $\beta_{h_{0}}(\rho_{0})=0$. 
According to the proof of Lemma \ref{lem:seceqn1}, we have 
\[\langle \,e\,|\,\phi_{F,h_{0}}(\rho)(s_{\xi}(\rho))\,\rangle 
=\sqrt{w_{z}}\langle\, \Omega_{\rho^{'}}^{h_{0}}\,|\,
\pi(\langle \xi^{'}|\xi\rangle )(z+h_{0})\,\rangle \]
for $z=\beta_{h_{0}}(\rho)$, $\rho \in {\cal V}_{h_{0}}$.
For an anti-holomorphic tangent vector $\bar{Y}$ of ${\cal P}_{b}$, 
we have 
\[\bar{\partial}_{\bar{Y}}\phi_{F,h}\left(s_{\xi}(\rho)\right)
={\cal O}([-2^{-1}w_{z}\langle Y|z\rangle \cdot\xi  ]_{\rho},\, h\,)\]
from which follows
$\left.\bar{\partial}_{\bar{Y}}\phi_{F,h}(\rho)\left(
s_{\xi}(\rho)\right)\right|_{z=0}=0$.
We see that the anti-holomorphic derivative of $s_{\xi}$
vanishes at each point in ${\cal P}_{b}$. 
Hence $s_{\xi}$ is holomorphic.

\noindent (ii)
For $z\in {\cal H}_{h}$, we have 
\[\{f_{A}\circ\beta_{h}^{-1}\}(z)=w_{z}\cdot \langle (z+h)|\pi(A)(z+h)\rangle.\]
Then the representation $X_{f_{A}}^{h}$ of the Hamiltonian vector field
$X_{f_{A}}$ of $f_{A}$ at $({\cal V}_{h},\beta_{h},{\cal H}_{h})$ is 
\[(X_{f_{A}}^{h})_{z}=-\sqrt{-1}\{
\pi(A)(z+h)-\langle h|\pi(A)(z+h)\rangle (z+h)\}\quad(z\in {\cal H}_{h}).\]
If we take $h$ such that $\beta_{h}(\rho_{0})=0$, then it holds that 
\[(X_{f_{A}}^{h})_{0}=-\sqrt{-1}\{ \pi(A)h-\langle h|\pi(A)h\rangle h\}.\]
The connection $D$ satisfies
$\langle\, v\,|\,(D_{X_{f_{A}}}s)(\rho_{0})\,\rangle _{\rho_{0}}
=\partial_{\rho_{0}}(
\langle v|s(\cdot)\rangle _{\rho_{0}})(X_{ f_{A} })$
for $v\in {\cal E}_{X,\rho_{0}}$ and $s\in \Gamma_{\infty}(\exx)$.
Hence  we have 
$(D_{X_{f_{A}}}s_{\xi})(\rho_{0})=[\, \xi a_{ X_{ f_{A} } ,0}\,]_{\rho_{0}}$
where $a_{X_{f_{A}},0}\in {\cal A}$ satisfies that
\[ \pi(a_{X_{f_{A}},0})h=X_{f_{A}}= 
-\sqrt{-1}(\,\pi(A)-\langle h|\pi(A)h\rangle \,)h.\]
Therefore we have 
$\sqrt{-1}(D_{X_{f_{A}}}s_{\xi})(\rho_{0})
= s_{\xi A}(\rho_{0})-s_{\xi}(\rho_{0})f_{A}(\rho_{0})$
from which follows 
\[(s_{\xi}*f_{A})(\rho_{0})=s_{\xi}(\rho_{0})f_{A}(\rho_{0})
+ \sqrt{-1}(D_{X_{f_{A}}}s_{\xi})(\rho_{0})= s_{\xi A}(\rho_{0}).\]
Therefore we obtain the statement.
\qedh

Finally, we come to prove Theorem \ref{Thm:main00}.
\\
\\
\noindent
{\it Proof of Theorem \ref{Thm:main00}.}
(i) 
By definition, we see that $\gax=\Psi(X)$.
Therefore the statement follows from Lemma \ref{lem:mlem} (i).  

\noindent
(ii) 
Because $\gax=\Psi(X)$, $\kup=f({\cal A})$ and Lemma \ref{lem:mlem} (ii)
for $D$, the linear space $\gax$ is a right $\kup$-module.

Because 
$\rho(\langle \xi|\xi^{'}\rangle )=f_{\langle \xi|\xi^{'}\rangle }(\rho)$,
we see that $H(\Psi(\xi),\, \Psi(\xi^{'}))
=f_{\langle \xi|\xi^{'}\rangle }\in {\cal K}_{u}({\cal P})$.
Hence  $H(s,s^{'})\in {\cal K}_{u}({\cal P})$
for each $s, s^{'}\in \Gamma_{X}$.
For $\xi,\eta\in X$ and $A\in {\cal A}$,
we can verify that $H_{\rho}(\, s_{\eta},\, \,s_{\xi}*f_{A}\, )
= \{H(s_{\eta}, \,s_{\xi })*f_{A}\}(\rho)$ 
where we use $H_{\rho}\left(\Psi(\xi),
\Psi(\eta)\right)=\rho(\langle \xi|\eta\rangle )$
for $\xi,\eta\in X$ and $\rho\in {\cal P}$.
Hence $H(s, \, s^{'}*l)=H(s, \, s^{'}) *l$ for each
$s,s^{'}\in \gax$ and $l\in \kup$.
From the property of the ${\cal A}$-valued inner
product of $X$ and by the proof of Lemma \ref{lem:mlem} (i), we obtain
$\|H(s,s)\|^{1/2}=\|s\|$ for each $s\in \gax$
where the norm of $H(s,s)$ is the one defined in (\ref{eqn:normtwo}).
Hence the statement holds.

\noindent
(iii)
Because $H(\Psi(\xi),\, \Psi(\xi^{'}))
=f_{\langle \xi|\xi^{'}\rangle }$,
the map $\Psi$ is an isometry from $X$ onto $ \gax$.
Rewrite module actions $\phi$ and $\psi$ on $X$ and $ \gax$, 
respectively, by
\[\phi(\xi,A)\equiv \xi A,\quad \psi(s, l)\equiv s*l
\quad(\xi\in X,\,A\in {\cal A},\,s\in \gax,\,l\in \kup).\]
Then we obtain that $\psi\circ (\Psi\times f)=\Psi \circ \phi$
by Lemma \ref{lem:mlem} (ii). 
Hence the statement holds.
\qedh

\noindent{\bf Acknowledgement:}
The author would like to thank Prof. Izumi Ojima
and Takeshi Nozawa for a critical reading of this paper.
We are also grateful to Prof. George A. Elliott for his helpful advice.

\vv

\noindent
{\Large {\bf Appendix}}

\appendix

%
%
\sftt{Example of uniform K\"{a}hler bundle}
%
%
\begin{ex}
\label{ex:exx}
{\rm
Assume that ${\cal H}$ is a separable infinite dimensional Hilbert space.
\begin{enumerate}
\item 
Let ${\cal A}\equiv{\cal L}({\cal H})$
be the C$^{*}$-algebra of all bounded linear operators on ${\cal H}$.
The uniform K\"{a}hler bundle of ${\cal A}$ is $({\cal P}({\cal H})
\cup{\cal P}_{-}, \,
p, \,
2^{[0,1]}\cup\{b_{0}\})$
where ${\cal P}({\cal H})$ is the projective Hilbert space of ${\cal H}$,
${\cal P}_{-}$ is the union of a family of projective Hilbert spaces indexed 
by the power set of the closed interval $[0,1]$
and $\{b_{0}\}$ is the one-point set corresponding to
the equivalence class of identity representation 
$({\cal H},id_{{\cal L}({\cal H})})$ of ${\cal L}({\cal H})$ 
on ${\cal H}$. Since the primitive spectrum of ${\cal L}({\cal H})$
is a two-point set, the topology of $2^{[0,1]}\cup\{b_{0}\}$ is equal to
$\{\, \emptyset,\,2^{[0,1]},\,\{b_{0}\},\,2^{[0,1]}\cup\{b_{0}\}\,\}$ \cite{K-R}.
In this way, the base space of the uniform K\"{a}hler bundle is
not always a singleton when the C$^{*}$-algebra is type $I$.
\item
For the C$^{*}$-algebra ${\cal A}$ generated by the Weyl form
of the $1$-dimensional canonical commutation relation
$U(s)V(t)=e^{\sqrt{-1}st}V(t)U(s)$ for $s,t\in {\bf R}$,
its uniform K\"{a}hler bundle is $({\cal P}({\cal H}),p,\{1pt\})$.
The spectrum is a one-point set $\{1pt\}$ from von Neumann
uniqueness theorem \cite{BraRobi}.
\item
The {\it CAR algebra} ${\cal A}$ is a UHF algebra with the nest
$\{M_{2^{n}}({\bf C})\}_{n\in {\bf N}}$.
The uniform K\"{a}hler bundle has the base space $2^{{\bf N}}$
and each fiber on $2^{{\bf N}}$ is a separable infinite
dimensional projective Hilbert space where $2^{{\bf N}}$
is the power set of the set ${\bf N}$ of all natural numbers
with trivial topology, that is, the topology of $2^{{\bf N}}$
is just $\{\emptyset,2^{{\bf N}}\}$.
In general, the Jacobson  topology of the spectrum of 
a simple C$^{*}$-algebra is trivial \cite{K-R}.
\end{enumerate}
}
\end{ex}

%

\end{document}